\documentclass[11pt]{article}
\usepackage{graphics,amsmath,amssymb}
\usepackage{latexsym}
\usepackage{epsfig}
\usepackage{hyperref}
\usepackage{fancybox}

\def\epsilon{\varepsilon}

\def\phi{\varphi}

\def\om{\omega}

\newtheorem{theorem}{Theorem}[section]
\newtheorem{lemma}[theorem]{Lemma}

\newtheorem{corollary}[theorem]{Corollary}
\newtheorem{definition}[theorem]{Definition}

\newtheorem{proposition}[theorem]{Proposition}
\newtheorem{remark}[theorem]{Remark}

\def\N{{\mathbb N}}
\def\C{{\mathbb C}}
\def\R{{\mathbb R}}

\newenvironment{Proof}{\removelastskip\par\medskip
\noindent{\em Proof.} \rm}{\penalty-20\null$\square$\par\medbreak}

\title{\bf A semilinear integro-differential 
equation: 
\\
global existence and hidden regularity}

\author{Paola Loreti
\thanks{Dipartimento di Scienze di Base e Applicate per l'Ingegneria,
Sapienza Universit\`a di Roma,
Via Antonio Scarpa 16, 00161 Roma (Italy); e-mail: 
$<$paola.loreti@sbai.uniroma1.it$>$ }
\and Daniela Sforza
\thanks{Dipartimento di Scienze di Base e Applicate per l'Ingegneria, 
Sapienza Universit\`a di Roma,
Via Antonio Scarpa 16, 00161 Roma (Italy); e-mail: 
$<$daniela.sforza@sbai.uniroma1.it$>$ }}

\begin{document}
\date{}

\maketitle

\begin{abstract}

Here  we  show a hidden regularity result for  nonlinear wave equations
with an integral term of convolution type  and Dirichlet boundary conditions. Under    general assumptions on the nonlinear term  and on the integral kernel we are able to state results about global  existence of strong and mild  solutions without any further smallness on the initial data. 
Then we  define the trace of the normal derivative of the solution showing a regularity result. In such a way we extend to integrodifferential equations with  nonlinear term well-known results available in the literature for linear wave equations with memory.

\end{abstract}

%{\bf Keywords:} 

\bigskip
\noindent
\section{Introduction}

Let $\Omega\subset\R^N$  $(N\ge1)$ be a bounded open domain of class $C^2$.
Let us denote by $\nu$ the outward unit normal vector to the boundary
$\Gamma$.
In this paper
we will consider the Cauchy problem for nonlinear wave equations
with a general integral term  and Dirichlet boundary conditions:
\begin{equation}\label{eq:cauchyI}
\begin{cases}
\displaystyle
 u_{tt}(t,x) =\triangle u (t,x)+\int_0^t \dot{a} (t-s)\triangle u(s,x)\,ds+g(u(t,x)),
\quad t\ge0,\,\, x\in \Omega,
\\
u(t,x)=0\qquad   t\ge0, \,\, x\in\Gamma,
\\
u(0,x)=u_{0}(x),\quad
u_t(0,x)=u_{1}(x),\qquad  x\in \Omega.
\end{cases}
\end{equation}
According to the physical model as proposed in \cite{RHN}, we will assume that the integral kernel satisfies
\begin{equation}\label{eq:intro-a}
\begin{split}
& a: (0,\infty)\to\R\ \text{ is  a positive definite function with}\  a (0)<1,
 \\
 & a\,, \dot{a}\in L^1(0,+\infty), 
  \end{split}
\end{equation}
and the nonlinear term fulfils the following conditions:
\begin{itemize}
\item
 $g\in C(\R)$ such that there exist $\alpha\ge0$, with $\alpha(N-2)\le2$, and $C>0$ so that
\begin{equation}\label{eq:intro-g}
\begin{split}
& g(0)=0,
\\
&|g(x)-g(y)|\le C(1+|x|^\alpha+|y|^\alpha)|x-y|
\qquad\forall x,y\in\R
\,,
\end{split}
\end{equation}
\item 
set $\displaystyle G(t)=\int_0^t g(s)\ ds$, there exists $C_0>0$
such that
\begin{equation}\label{eqn:hypnn0}
G(t) \le C_0|t|^2
\qquad\forall t\in\R\,.
\end{equation} 
\end{itemize}
%As well known, for more regular initial data a  mild solution is a strong solution and
% in the non linear case the notions of mild solution and weak solution are not equivalent.
% We will show a global existence result about mild and strong solutions.
%
 We will establish the following global existence result without any smallness assumption on initial data.
\begin{theorem}
Under the assumptions
\eqref{eq:intro-a}--\eqref{eqn:hypnn0},
for any  $(u_0,u_1)\in H^1_0(\Omega)\times L^2(\Omega)$ 
problem $\eqref{eq:cauchyI}$ admits a
unique  mild solution $u$ on $[0,\infty )$. 
\end{theorem}

In our previous work \cite{LoretiSforza2}  we study  the linear case of \eqref{eq:cauchyI} where
the integral kernel 
$\dot{a}: [0,\infty)\to (-\infty,0]$ is a locally absolutely continuous
function,
$\dot{a}(0)<0$,
$\ddot{a}(t)\ge 0$ for a.e. $ t\ge 0$ and
$a(0)<1$.
%
%
%is given by $\dot{a} (t)=-k(t)$ where
%$k: [0,\infty)\to [0,\infty)$ is a locally absolutely continuous
%function,
%$k(0)> 0$,
%$\dot{k}(t)\le 0$ for a.e. $ t\ge 0$ and
%$\int_0^{\infty} k (t)\ dt<1$.

In this paper  the existence result   may  be stated for more general kernels, as

\begin{itemize}
\item 

$a(t)=a_0\int_t^\infty \frac{e^{-\alpha s}}{s^\beta}\ ds$, with $\alpha>0$, $0\le\beta<1$ and $0\le a_0<\frac{\Gamma(1-\beta)}{\alpha^{1-\beta}}$\,,
\item 

$a(t)=\int_t^\infty (a_0 s+a_1)e^{-\alpha s}\ ds=
(\frac{a_0}{\alpha}t+\frac{a_0+\alpha a_1}{\alpha^2})e^{-\alpha t}$, 

with $\alpha>0$, $a_0,a_1\ge0$, $\frac{a_0+\alpha a_1}{\alpha^2}<1$, $\alpha a_1-a_0\ge0$,

\item
$a(t)=k\int_t^\infty  \frac 1 {(1+s)^\alpha}ds,\,\,\, $ with $ k>0$ such that $a(0)<1, \,\, \alpha>2.$
\end{itemize}
Examples of fading memory kernels can be found in  \cite{DE},\cite{Ka}.

Some examples of  $g$ satisfying assumptions \eqref{eq:intro-g}--\eqref{eqn:hypnn0} are
$$
g(x)=c \vert x\vert^p x,\quad c<0,\ \  p(N-2)\leq 2\,, \qquad\qquad g(x)=c\sin x,\quad c\in\R\,. 
$$
%We will define 
%the energy of a mild solution $u$  in a given interval $[0,T]$ as follows
%\begin{equation*}
%E(t)=\;\frac{1}{2}\int_{\Omega}|u_t|^2\ dx+\frac{1-a(0)}{2}\int_{\Omega} |\nabla u|^2\ dx-\int_{\Omega} G(u)\ dx
%\,.
%\end{equation*}
In addition, for more regular kernels we will  prove a so-called hidden regularity result.
\begin{theorem}
Assume
\eqref{eq:intro-g}--\eqref{eqn:hypnn0},
\begin{equation*}
C_0<\lambda(1-a(0))/2
\,,
\qquad
\lambda:=\inf\{\|\nabla v\|_{L^2}^2,\ v\in H^1_0(\Omega),\ \| v\|_{L^2}=1\},
\end{equation*}
and
\begin{equation}\label{eq:lessr}
\begin{split}
&a\in C^1 ([0,\infty)),\ \dot{a}(0)<0,\ a(t)\ge0,\ \dot{a}(t)\le0\ \forall t\ge0,
\\
& \ddot{a}(t)\in L^1_{loc}(0,+\infty), \ddot{a}(t)\ge0, \ \text{a.e.}\ t\ge0,
\\
&a\,, \dot{a}\in L^1(0,+\infty),\ a(0)<1\,.
\end{split}
\end{equation} 
Let $T>0$, there exists a constant $c=c(T)>0$ such that for any $u_0\in
H^1_0(\Omega)$ and $u_1\in L^{2}(\Omega)$ if $u$ is the mild solution  of
\eqref{eq:cauchyI}, then, denoting by  $\partial_\nu u$ the normal derivative, we have
\begin{equation*}
\int _{0}^T \int _\Gamma |\partial_\nu u|^2\ d\Gamma \
dt\le c(\|\nabla u_0\|^2_{L^2}+\|u_1\|^2_{L^2})
\,.
\end{equation*}
Moreover, if the energy $E(t)$ of the solution $u$ satisfies 
\begin{equation*}
\int_0^t E(s)\ ds\leq c_0 E(0)\qquad\forall t \geq 0
\qquad
(c_0>0\ \text{ independent of  }\ t),
\end{equation*}
then we have
\begin{equation*}
\partial_\nu u\in L^2(0,\infty;L^2(\Gamma )).
\end{equation*}
  \end{theorem}
  
 The proof of the existence of the solution $u$ relies on energy estimates.
 Although we use some  results  obtained in \cite {ACS}, \cite {CS2} and \cite {CS3},
 here we are interested to treat initial data without any smallness, so the previous results have to be adapted     
 to consider our  setting. We  also mention the papers \cite{Berrimi,Kim}.

To understand how to frame our paper in the literature,
we recall briefly some known results.
Lasiecka and Triggiani  \cite {LasTri}  established  the  hidden regularity
property for the weak solution $u$ of the wave equation with Dirichlet boundary conditions
%\begin{equation*}
%\begin{cases}
%\displaystyle
% u_{tt}(t,x) =\triangle u (t,x),
%\quad t\ge0,\,\, x\in \Omega,
%\\
%u(t,x)=0\qquad   t\ge0, \,\, x\in\Gamma,
%\\
%u(0,x)=u_{0}(x),\quad
%u_t(0,x)=u_{1}(x),\qquad  x\in \Omega,
%\end{cases}
%\end{equation*}
that is 
\begin{equation*}
\partial_\nu u\in L^2_{loc}(\R;L^2(\Gamma )).
\end{equation*}
The  term hidden was proposed  by J.L. Lions  \cite {Lio} for the wave  equation
 in the context of the exact controllability problems.
 Later in  \cite {Lio2} J.L. Lions proved that the weak solution of the nonlinear wave  equation 
 \begin{equation*}
 u_{tt}(t,x) =\triangle u (t,x) -\vert u\vert^p u,
\quad t\ge0,\,\, x\in \Omega,
\end{equation*}
%\begin{equation*}
%\begin{cases}
%\displaystyle
% u_{tt}(t,x) =\triangle u (t,x) -\vert u\vert^p u,
%\quad t\ge0,\,\, x\in \Omega,
%\\
%u(t,x)=0\qquad   t\ge0, \,\, x\in\Gamma,
%\\
%u(0,x)=u_{0}(x),\quad
%u_t(0,x)=u_{1}(x),\qquad  x\in \Omega,
%\end{cases}
%\end{equation*}
satisfies a trace regularity result. Milla Miranda and Medeiros \cite {MM}
enlarged the class of nonlinear terms by means of approximation arguments. 
However they do not consider memory terms in the equation, that is $\dot a\equiv0$. To our knowledge it seems that there are not previous papers studying the hidden regularity for solutions of nonlinear integro-differential problems when the integral kernels satisfy the assumptions \eqref{eq:lessr}.

The plan of our paper is the following. In Section 2 we list some notations and preliminary results. In Section 3 we establish existence and uniqueness results of mild and strong solutions. Finally, in Section 4 we give hidden regularity results for a nonlinear equation with memory.

\section{Preliminaries}
Let 
$
L^2(\Omega)
$
be endowed with the usual inner product and norm 
$$
\|u\|_{L^2}=\left(\int_{\Omega}|u(x)|^{2}\ dx\right)^{1/2}\qquad
u\in L^2(\Omega)\,.
$$
Throughout the paper we will  use a standard notation for the integral convolution between two functions, that is 
\begin{equation}
h*u(t):=\int_0^th(t-s) u(s)\ ds
\,.
\end{equation}
A well-known result concerning integral equations (see e.g. \cite[Theorem 2.3.5]{GLS}), that we will use later is the following.
\begin{lemma} \label{le:unicita}
Let $h\in L^1 (0,T)$, $T>0$. If the function $\varphi(t)+h*\varphi(t)$ belongs to $L^2 (0,T;L^2(\Omega))$
then $\varphi\in L^2 (0,T;L^2(\Omega))$  and
there exist a positive constant $c_1=c_1(\|h\|_{L^1(0,T)})$, depending on the norm $\|h\|_{L^1(0,T)}$, such that
\begin{equation}\label{eq:unicita}
\int_0^T\|\varphi(t)\|_{L^2}^2\ dt\le
c_1\int_0^T\big\|\varphi(t)+h*\varphi(t)\big\|_{L^2}^2\ dt
%\le c_2\int_0^T\|\varphi(t)\|_{L^2}^2\ dt
\,.
\end{equation}
\end{lemma}
Recall that $h$ is a {\em positive definite kernel}  if for any $y\in L^2_{loc}(0,\infty;L^2(\Omega))$ we have
\begin{equation}\label{postype}
\int_0^t\int_\Omega y(\tau,x)\int_0^\tau h(\tau-s)y(s,x)\ ds \ dx\ d\tau\ge 0\,,
\qquad t\ge 0\,.
\end{equation}
Also, $h$ is said to be a {\em strongly positive definite kernel} if there exists a constant $\delta>0$ such
that 
$h(t)-\delta e^{-t}$ is  positive definite. This stronger notion for the integral kernel allows to obtain uniform estimates for solutions of integral equations, see \cite[Corollary 2.12]{CS3}. For completeness we recall here that result, because we will use it later.
\begin{lemma} \label{le:unicita2}
Let $a\in L^1 (0,\infty)$ be a strongly positive definite kernel such that  $\dot a\in L^1 (0,\infty)$ and $a(0)<1$. If the function $\varphi(t)+\dot a*\varphi(t)$ belongs to $L^2 (0,\infty;L^2(\Omega))$
then $\varphi\in L^2 (0,\infty;L^2(\Omega))$ 
and there exist a positive constant $c_1$, such that
\begin{equation}\label{eq:unicita2}
\int_0^\infty\|\varphi(t)\|_{L^2}^2\ dt\le
c_1\int_0^\infty\big\|\varphi(t)+\dot a*\varphi(t)\big\|_{L^2}^2\ dt
%\le c_2\int_0^\infty\|\varphi(t)\|_{L^2}^2\ dt
\,.
\end{equation}
\end{lemma}

Regarding the nonlinear term, we will follow the approach pursued in \cite{CH} for the nonintegral case when $\dot a\equiv0$. Precisely, we will consider a function $g\in C(\R)$ such that there exist $\alpha\ge0$, with $(N-2)\alpha\le2$, and $C>0$ so that
\begin{equation}\label{eq:g1}
\begin{split}
& g(0)=0,
\\
&|g(x)-g(y)|\le C(1+|x|^\alpha+|y|^\alpha)|x-y|
\qquad\forall x,y\in\R
\,.
\end{split}
\end{equation}
In \cite[Proposition 6.1.5]{CH} the following result has been proved.
\begin{proposition}\label{eq:lipcon}
If $g$ satisfies the hypotheses \eqref{eq:g1}, then $g$ is Lipschitz continuous from bounded subsets of $H_0^1(\Omega )$ to $L^2(\Omega)$.
In particular, there exists a positive constant  $C$ such that 
\begin{equation}\label{eq:lipcon1}
\int_\Omega |g(u(x))|^2\ dx\le C\int_\Omega |\nabla u(x)|^2\ dx
\qquad
\forall u\in H_0^1(\Omega )
\,.
\end{equation}

\end{proposition}
We will assume that the integral kernel satisfies the following conditions:
\begin{equation}\label{eq:k1}
\begin{split}
&a: (0,\infty)\to\R \quad\text{is  a positive definite function,}
\\
&a\,, \dot{a}\in L^1(0,+\infty),
\\
& a(0)<1\,.
\end{split}
\end{equation}
For reader's convenience we begin with recalling some known notions and results. 
First, we  write the Laplacian as an abstract operator.
Indeed, we define the operator $A:D(A)\subset L^2(\Omega)\to L^2(\Omega)$ as 
$$
\begin{array}{l}
D(A)=H^2(\Omega)\cap H_0^1(\Omega) \\
\\
Au(x)=\displaystyle -\Delta u(x)\qquad u\in D(A)\,,\; x\in \Omega \;\mbox{a.e.}
\end{array}
$$
We recall, see e.g. \cite[Definition 3.1]{CS2}, that there exists a unique family $\{{\cal R}(t)\}_{t\ge 0}$ of bounded linear operators in $L^2(\Omega)$ the so-called 
resolvent for the linear equation
\begin{equation}\label{eq:linear}
u''(t)+Au (t)+\int_0^t\dot a (t-s)Au(s)\,ds=0\,,
\end{equation}
that satisfy
the following conditions:
\begin{description}
\item[\rm(i)]
${\cal R}(0)$ is the identity operator and ${\cal R}(t)$ is strongly continuous on $[0,\infty)$,
that is, for all $u\in L^2(\Omega)$, ${\cal R}(\cdot)u$ is continuous;
\item[\rm(ii)]
${\cal R}(t)$ commutes with $A$, which means that ${\cal R}(t)D(A)\subset D(A)$ and
$$
A{\cal R}(t)u={\cal R}(t)Au\,,\qquad u\in D(A)\,,\,t\ge 0\,;
$$
 \item[\rm(iii)]
for any  $u\in D(A)$, ${\cal R}(\cdot)u$ is twice continuously
differentiable in $L^2(\Omega)$ on $[0,\infty)$ and $ {\cal R}'(0)u=0$;
\item[\rm(iv)]
for any $u\in D(A)$ and any $t\ge 0$, 
\begin{equation*}
 {\cal R}''(t)u+A{\cal R}(t)u+\int_0^t\dot a(t-\tau)A{\cal R}(\tau)u\,d\tau=0\,.
\end{equation*}
\end{description}
In the sequel we will use the following uniform estimates for the
resolvent, see e.g. \cite[Proposition 3.4-(i)]{CS2}, taking into account that $
D(A^{1/2}) = H^1_0(\Omega)
$.
\begin{proposition}
For any $u\in L^2(\Omega)$ and any $t> 0$, we have $1*{\cal R}(t)u\in H_0^1(\Omega)$ and
\begin{eqnarray}\label{stima1}
\|{\cal R}(t)u\|_{L^2}^2+
\big(1-a(0)\big)
\left\|\nabla(1*{\cal R})(t)u\right\|_{L^2}^2
%_{{\mathcal L}(X)}
\le \|u\|_{L^2}^2\,.
\end{eqnarray}
In particular, $\nabla(1*{\cal R})(\cdot)$ is strongly continuous in $L^2(\Omega)$.
\end{proposition}
Let $0<T\le\infty$ be given. 
We recall some notions of solution for the semilinear equation
\begin{equation}\label{eq:stato0}
 u_{tt}(t,x) =\triangle u (t,x)+\int_0^t \dot{a} (t-s)\triangle u(s,x)\,ds+g(u(t,x)),
\quad t\in[0,T]\,,\, x\in \Omega\,.
\end{equation}
\begin{definition}
 We say that $u$ is a {\bf strong solution} of \eqref{eq:stato0} on $[0,T]$ if
$$u\in C^2([0,T];L^2(\Omega))\cap C([0,T];H^2(\Omega)\cap H^1_0(\Omega))$$ and $u$ satisfies 
\eqref{eq:stato0} for every $t\in[0,T]$.

Let $u_0$, $u_1\in L^2(\Omega)$.  A function 
$u\in  C^1([0,T];L^2(\Omega))\cap C([0,T];H^1_0(\Omega))$ 
is a {\bf mild solution} of  \eqref{eq:stato0} on $[0,T]$ with initial conditions
\begin{equation}\label{eq:inicond}
u(0)=u_{0},\quad
u_t(0)=u_{1},
\end{equation}
if
\begin{equation}\label{eq:varparn}
u(t)={\cal R}(t)u_0+\int_0^t{\cal R}(\tau)u_1d\tau+\int_0^t1*{\cal R}(t-\tau)g(u(\tau))d\tau,
\end{equation}
where $\{{\cal R}(t)\}$ is the resolvent for the linear equation \eqref{eq:linear}.
\end{definition}
Notice that the convolution term in (\ref{eq:varparn}) is well defined, thanks
to Proposition \ref{eq:lipcon}.
A strong solution is also a mild one.

Another useful notion of generalized solution of \eqref{eq:stato0} is the so-called {\em weak solution}, that is a function 
$u\in C^1([0,T];L^2(\Omega))\cap C([0,T];H^1_0(\Omega))$ such that for any $v\in H^1_0(\Omega)$, $ t\to\int_\Omega u_t v\ dx\in C^1([0,T])$ 
and 
\begin{equation}\label{eq:weak_sol}
\frac {d}{dt} \int_\Omega u_t v\ dx 
=-\int_\Omega \nabla u\cdot\nabla v\ dx
-\int_\Omega\int_0^t \dot a(t-s)\nabla u(s)\,ds\cdot\nabla v\ dx
+\int_\Omega g(u(t))v\ dx\,,
\qquad\forall t\in[0,T]\,.
\end{equation}
Adapting a classical argument due to Ball \cite{Ball}, one can show that any mild solution of \eqref{eq:stato0}
is also a weak solution, and the two notions of solution are equivalent in the linear case when $g\equiv 0$ (see also \cite{Pruss}).

Throughout the paper we denote with the symbol $\cdot$ the Euclidean scalar product in $\R^N$.

\section{Existence and uniqueness of mild and strong solutions}

The next proposition ensures the local existence and uniqueness
of the mild solution for the Cauchy problem
\begin{equation}\label{eq:cauchy2}
\begin{cases}
\displaystyle
 u_{tt}(t,x) =\triangle u (t,x)+\int_0^t \dot{a} (t-s)\triangle u(s,x)\,ds+g(u(t,x)),
\quad t\ge0,\,\, x\in \Omega,
\\
u(t,x)=0\qquad   t\ge0, \,\, x\in\Gamma,
\\
u(0,x)=u_{0}(x),\quad
u_t(0,x)=u_{1}(x),\qquad  x\in \Omega.
\end{cases}
\end{equation}
The proof relies on  suitable regularity estimates for the resolvent $\{{\cal R}(t)\}$ such as \eqref{stima1}
%and well-posedness for the linear problem 
(for more details see e.g.
\cite[section 3]{CS2}) and a standard fixed point argument (see
\cite{CS} for an analogous proof).

\begin{proposition}\label{cir}
If
$u_0\in H^1_0(\Omega)$ and $u_1\in L^2(\Omega)$,
there exists a positive number
$T$ such that the Cauchy problem \eqref{eq:cauchy2}
admits a
unique mild solution on $[0,T]$. 

\end{proposition}

Assuming more regular data and using standard argumentations, one can show that the mild solution is a strong one.

\begin{proposition}\label{th:t2}
Let  $u_0\in H^2(\Omega)\cap H^1_0(\Omega)$ and $u_1\in H^1_0(\Omega)$. Then, the mild
solution of the Cauchy problem \eqref{eq:cauchy2} in $[0,T]$
%$T_0\in (0,T]$, 
is a strong solution.
In addition, $u$ belongs to
$C^1([0,T];H^1_0(\Omega))$.
\end{proposition}

To investigate the existence for all $t\ge0$ of
the solutions,
for $g$ satisfying \eqref{eq:g1} we introduce $G\in C(\R)$ by means of
\begin{equation}
G(t)=\int_0^t g(s)\ ds
\,.
\end{equation}
We define the energy of a mild solution $u$ of \eqref{eq:cauchy2} on
a given interval $[0,T]$, as
\begin{equation}\label{eqn:energy}
E(t)=\;\frac{1}{2}\int_{\Omega}|u_t|^2\ dx+\frac{1-a(0)}{2}\int_{\Omega} |\nabla u|^2\ dx-\int_{\Omega} G(u)\ dx
\,.
\end{equation}
In view of \eqref{eq:g1} we have
\begin{equation}\label{eq:Gu-0}
\int_\Omega |G(u_0(x))|\ dx\le C\int_\Omega|\nabla u_0(x)|^2\ dx
\qquad
\forall u_0\in H^1_0(\Omega)
\,,
\end{equation}
and hence
\begin{equation}
E(0)\le C(\|\nabla u_0\|_{L^2}^2+\| u_1\|_{L^2}^2)
\qquad
u_0\in H^1_0(\Omega),\ u_1\in L^2(\Omega)
\,.
\end{equation}
About the energy of the solutions, we recall some known results, see \cite[Lemma 3.5]{CS3}.
\begin{lemma}\label{le:diss}
\begin{itemize}
\item [i)]
If $u_0\in H^2(\Omega)\cap H^1_0(\Omega)$ and $u_1\in H^1_0(\Omega)$, then the strong solution $u$  of problem 
$\eqref{eq:cauchy2}$ on $[0,T]$
satisfies the identity 
\begin{multline}\label{eq:disss}
E(t)+\int_0^t\int_\Omega a*\nabla u_t(s)\cdot \nabla u_t(s)\ dx\ ds
\\
= E(0)+a(0)\int_\Omega |\nabla u_{0}|^2\ dx
-a(t)\int_\Omega \nabla u_{0}\cdot \nabla u(t)\ dx
-\int_0^t\dot a(s)\int_\Omega \nabla u_{0}\cdot \nabla u(s)\ dx\ ds
\,,
\end{multline}
for any
$t\in [0,T]$.
\item [ii)]
If
$u_0\in H^1_0(\Omega)$ and $u_1\in L^2(\Omega)$, then the mild solution $u$  of problem $\eqref{eq:cauchy2}$ on $[0,T]$ verifies 
\begin{equation}\label{eq:diss}
E(t)
\le E(0)+a(0)\int_\Omega |\nabla u_{0}|^2\ dx
-a(t)\int_\Omega \nabla u_{0}\cdot \nabla u(t)\ dx
-\int_0^t\dot a(s)\int_\Omega \nabla u_{0}\cdot \nabla u(s)\ dx\ ds
\,,
\end{equation}
for any
$t\in [0,T]$.
\end{itemize}
\end{lemma}
Assuming an extra condition on  $G$, global existence will follow 
for all data. For further convenience we introduce the notation
\begin{equation}\label{eq:lambda}
\lambda=\inf\{\|\nabla v\|_{L^2}^2,\ v\in H^1_0(\Omega),\ \| v\|_{L^2}=1\}
\,.
\end{equation}

\begin{theorem}\label{th:alldata}
Suppose that there exists $C_0>0$
such that 
\begin{equation}\label{eqn:hypnn}
G(t) \le C_0|t|^2
\qquad\forall t\in\R\,.
\end{equation}
Then  for any  $(u_0,u_1)\in H^1_0(\Omega)\times L^2(\Omega)$ 
problem $\eqref{eq:cauchy2}$ admits a
unique  mild solution $u$ on $[0,\infty )$. 

Moreover, if we suppose that the constant $C_0>0$ in \eqref{eqn:hypnn} satisfies 
\begin{equation}\label{eq:c0}
C_0<\lambda(1-a(0))/2,
\end{equation}
where $\lambda$ is defined in \eqref{eq:lambda}, then $E(t)$ is positive and  we have for any $t\ge 0$
\begin{eqnarray}
\label{eq:coercive}
& \displaystyle E(t)\ge \;\frac12\| u_t(t)\|_{L^2}^2+\frac{C}2\|\nabla u(t)\|_{L^2}^2\,, &
\\
\label{eq:boundE}
& 
\displaystyle
E(t)
\le C\big(\| u_{1}\|_{L^2}^2 +\| \nabla u_{0}\|_{L^2}^2\big)\,,&
\\
\label{eq:boundA}
& \| u_t(t)\|_{L^2}^2+\|\nabla u(t)\|_{L^2}^2
\displaystyle
\le C\big(\| u_{1}\|_{L^2}^2 +\| \nabla u_{0}\|_{L^2}^2\big)\,,&
\end{eqnarray}
%for any $t\ge 0$, where $ C(R)$
%is a positive, increasing, upper semicontinuous function such that $C(0)=0$.
where the symbol $C$ denotes positive constants, that can be different.

Furthermore, if $u_0\in H^2(\Omega)\cap H^1_0(\Omega)$ and $u_1\in H^1_0(\Omega)$, 
then $u$ is a strong solution of $\eqref{eq:cauchy2}$ on $[0,\infty)$, $u\in C^1([0,\infty);H^1_0(\Omega))$ and for any $t\ge 0$
\begin{equation}\label{eq:boundG*A}
E(t)+\int_0^t\int_\Omega a*\nabla u_t(s)\cdot\nabla u_t(s)\ dx\ ds
\le C\big(\| u_{1}\|_{L^2}^2 +\| \nabla u_{0}\|_{L^2}^2\big)\,.
\end{equation}
\end{theorem}
\begin{Proof}
Let $[0,T)$ be the maximal domain of the mild solution $u$ of $\eqref{eq:cauchy2}$. To prove  $T=\infty$, we will 
argue by contradiction and assume that $T$ is a positive real number. We will 
show that there exists a constant $C=C(T)>0$ such that
\begin{equation}\label{eq:energy}
\int_{\Omega}|u_t|^2\ dx+\int_{\Omega} |\nabla u|^2\ dx
\le C
\qquad\forall t\in[0,T)
\,.
\end{equation}
First, thanks to \eqref{eq:diss} we have
\begin{multline}\label{eq:bound1}
\int_{\Omega}|u_t|^2\ dx+\big(1-a(0)\big)\int_{\Omega} |\nabla u|^2\ dx
\\
\le \|u_1\|_{L^2}^2+\big(1-a(0)\big)\|\nabla u_0\|_{L^2}^2-2\int_{\Omega} G(u_0)\ dx
+2a(0)\|\nabla u_0\|_{L^2}^2
\\
-2a(t)\int_\Omega \nabla u_{0}\cdot \nabla u(t)\ dx
-2\int_0^t\dot a(s)\int_\Omega \nabla u_{0}\cdot \nabla u(s)\ dx\ ds
+2\int_{\Omega} G(u(t))\ dx
\\
\le \|u_1\|_{L^2}^2+\big(1+a(0)\big)\|\nabla u_0\|_{L^2}^2+2\int_{\Omega} |G(u_0)|\ dx
\\
-2a(t)\int_\Omega \nabla u_{0}\cdot \nabla u(t)\ dx
-2\int_0^t\dot a(s)\int_\Omega \nabla u_{0}\cdot \nabla u(s)\ dx\ ds
+2\int_{\Omega} G(u(t))\ dx
\,.
\end{multline}
We note that
\begin{multline*}
-2a(t)\int_\Omega \nabla u_{0}\cdot \nabla u(t)\ dx
\\
\le
2\|a\|_\infty\int_\Omega |\nabla u_{0}|\ |\nabla u(t)|\ dx
\le
\frac{1-a(0)}2\int_{\Omega} |\nabla u|^2\ dx
+\frac{2\|a\|_\infty^2}{1-a(0)}\|\nabla u_0\|_{L^2}^2
\,.
\end{multline*}
Putting the above estimate into \eqref{eq:bound1}, we obtain
\begin{multline}\label{eq:bound2}
\int_{\Omega}|u_t|^2\ dx+\frac{1-a(0)}2\int_{\Omega} |\nabla u|^2\ dx
\\
\le \|u_1\|_{L^2}^2+\Big(1+a(0)+\frac{2\|a\|_\infty^2}{1-a(0)}\Big)\|\nabla u_0\|_{L^2}^2+2\int_{\Omega} |G(u_0)|\ dx
\\
-2\int_0^t\dot a(s)\int_\Omega \nabla u_{0}\cdot \nabla u(s)\ dx\ ds
+2\int_{\Omega} G(u(t))\ dx
\,.
\end{multline}
Now, we have to estimate the last two terms on the right-hand side of the previous inequality. As regards the first one, we note that
\begin{equation}\label{eq:bound22}
-2\int_0^t\dot a(s)\int_\Omega \nabla u_{0}\cdot \nabla u(s)\ dx\ ds
\le
\frac{2\|\dot a\|_1}{1-a(0)}\|\nabla u_0\|_{L^2}^2
+\frac{1-a(0)}2\int_0^t|\dot a(s)|\int_\Omega |\nabla u(s)|^2\ dx\ ds
\,.
\end{equation}
Concerning the other integral,
assumption \eqref{eqn:hypnn}  yields for any $t\in [0,T)$
\begin{equation}\label{eqn:hypnn1} 
\int_\Omega G(u(t))\ dx\le C_0\int_\Omega |u(t)|^2\ dx
\,.
\end{equation}
In addition, we observe that
\begin{equation*}
\|u(t)\|_{L^2}^2=\|u_0\|_{L^2}^2+\int_0^t\frac{d}{ds}\int_\Omega |u(s)|^2\ dx\ ds
=\|u_0\|_{L^2}^2+2\int_0^t\int_\Omega u(s)u_t(s)\ dx\ ds
\,.
\end{equation*}
Since, by the definition \eqref{eq:lambda} of $\lambda$ we have 
\begin{equation}\label{eq:lambda2}
\int_\Omega |\nabla u|^2\ dx\ge \lambda\int_\Omega |u|^2\ dx
\,,
\end{equation}
we deduce
\begin{equation*}
\|u(t)\|_{L^2}^2\le
\|u_0\|_{L^2}^2+\int_0^t\Big(\int_\Omega |u_t(s)|^2\ dx+\frac1\lambda\int_\Omega |\nabla u(s)|^2\ dx\Big)\ ds
\,.
\end{equation*}
Therefore, by \eqref{eqn:hypnn1} 
\begin{multline}\label{eq:bound23}
\int_\Omega G(u(t))\ dx
\le C_0\|u_0\|_{L^2}^2
+C_0\int_0^t\Big(\int_\Omega |u_t(s)|^2\ dx+\frac1\lambda\int_\Omega |\nabla u_t(s)|^2\ dx\Big)\ ds
\\
\le C_0\|u_0\|_{L^2}^2
+M\int_0^t\Big(\int_\Omega |u_t(s)|^2\ dx+\frac{1-a(0)}2\int_\Omega |\nabla u(s)|^2\ dx\Big)\ ds
\,,
\end{multline}
where $M=C_0\max\{1, \frac{2\lambda}{1-a(0)}\}$.
Plugging \eqref{eq:bound22} and \eqref{eq:bound23} into \eqref{eq:bound2}, thanks also to \eqref{eq:Gu-0} we get
\begin{multline}\label{eq:bound3}
\int_{\Omega}|u_t|^2\ dx+\frac{1-a(0)}2\int_{\Omega} |\nabla u|^2\ dx
\\
\le \|u_1\|_{L^2}^2+\Big(1+a(0)+2\frac{\|a\|_\infty^2+\|\dot a\|_1}{1-a(0)}+\frac{C_0}{\lambda}+C\Big)\|\nabla u_0\|_{L^2}^2
\\
+\int_0^t\big(|\dot a(s)|+M\big)\Big(\int_\Omega |u_t(s)|^2\ dx+\frac{1-a(0)}2\int_\Omega |\nabla u(s)|^2\ dx\Big)\ ds\,.
\end{multline}
Applying Gronwall lemma, we obtain for any $t\in[0,T)$
\begin{multline*}
\int_{\Omega}|u_t|^2\ dx+\frac{1-a(0)}2\int_{\Omega} |\nabla u|^2\ dx
\\
\le \left(\|u_1\|_{L^2}^2+\Big(1+a(0)+2\frac{\|a\|_\infty^2+\|\dot a\|_1}{1-a(0)}+\frac{C_0}{\lambda}+C\Big)\|\nabla u_0\|_{L^2}^2\right)
e^{\int_0^t\big(|\dot a(s)|+M\big)\ ds}
\\
\le \left(\|u_1\|_{L^2}^2+\Big(1+a(0)+2\frac{\|a\|_\infty^2+\|\dot a\|_1}{1-a(0)}+\frac{C_0}{\lambda}+C\Big)\|\nabla u_0\|_{L^2}^2\right)
e^{\|\dot a\|_1+M T}
\,,
\end{multline*}
and hence, set 
\begin{equation*}
C(T)=\frac{e^{\|\dot a\|_1+M T}}{\min\{1,\frac{1-a(0)}2\}}\left(\|u_1\|_{L^2}^2+\Big(1+a(0)+2\frac{\|a\|_\infty^2+\|\dot a\|_1}{1-a(0)}+\frac{C_0}{\lambda}+C\Big)\|\nabla u_0\|_{L^2}^2
\right)
\end{equation*}
we have that \eqref{eq:energy} holds true.

To have a contradiction, we will prove that  $u\in  C([0,T];H^1_0(\Omega))\cap C^1([0,T];L^2(\Omega))$. 
First, set
\begin{equation}\label{eq:varparn10}
v(t)={\cal R}(t)u_0+\int_0^t{\cal R}(\tau)u_1d\tau
\qquad t\ge0
\,,
\end{equation}
we note that, thanks to the properties of the resolvent  we have
\begin{equation}\label{eq:v}
v(t)\in  C([0,+\infty);H^1_0(\Omega))\cap C^1([0,+\infty);L^2(\Omega))
\,.
\end{equation}
Since by \eqref{eq:varparn} and \eqref{eq:varparn10} we can write
\begin{equation}\label{eq:varparn11}
u(t)=v(t)+\int_0^t1*{\cal R}(t-\tau)g(u(\tau))d\tau,
\end{equation}
for $h>0$ and $0\le t<t+h<T$ we have
\begin{equation*}
\begin{split}
u(t+h)- u(t)=&v(t+h)-v(t)+\int_0^{t+h}1*{\cal R}(\tau)g(u(t+h-\tau))d\tau-\int_0^{t}1*{\cal R}(\tau)g(u(t-\tau))d\tau
\\
= &v(t+h)- v(t)+\int_0^{t}1*{\cal R}(\tau)\big[g(u(t+h-\tau))-g(u(t-\tau))\big]d\tau
\\
&+\int_t^{t+h}1*{\cal R}(\tau)g(u(t+h-\tau))d\tau
.
\end{split}
\end{equation*}
As a consequence, by \eqref{stima1} we have
\begin{equation*}
\begin{split}
\|\nabla u(t+h)- \nabla u(t)\|_{L^2}\le&\|\nabla v(t+h)-\nabla v(t)\|_{L^2}+\frac1{1-a(0)}\int_0^{t}\|g(u(s+h))-g(u(s))\|_{L^2}ds
\\
&+\frac1{1-a(0)}\int_0^{h}\|g(u(s))\|_{L^2}ds
.
\end{split}
\end{equation*}
Thanks to Proposition \ref{eq:lipcon} and \eqref{eq:energy} we deduce that
\begin{equation*}
\begin{split}
\|\nabla u(t+h)- \nabla u(t)\|_{L^2}\le&\|\nabla v(t+h)-\nabla v(t)\|_{L^2}+C\int_0^{t}\|\nabla u(s+h))-\nabla u(s)\|_{L^2}ds
\\
&+C\int_0^{h}\|\nabla u(s)\|_{L^2}ds
\\
\le&\|\nabla v(t+h)-\nabla v(t)\|_{L^2}+Ch+C\int_0^{t}\|\nabla u(s+h))-\nabla u(s)\|_{L^2}ds
,
\end{split}
\end{equation*}
where $C=C(T)>0$ is a positive constant. Applying Gronwall lemma, we get
\begin{equation*}
\begin{split}
\|\nabla u(t+h)- \nabla u(t)\|_{L^2}\le
\big(\|\nabla v(t+h)-\nabla v(t)\|_{L^2}+Ch\big)e^{C T}
,
\end{split}
\end{equation*}
and hence the function $\nabla u(t)$ is uniformly continuous in $[0,T[$ with values in $L^2(\Omega)$. Therefore $ u(t)$ can be also defined  in $T$ in a way that $u\in C([0,T];H^1_0(\Omega))$.
Moreover,  again by \eqref{eq:varparn11} we have
%for $h>0$ and $0\le t<t+h<T$ we have
%\begin{equation*}
%\begin{split}
%u_t(t+h)- u_t(t)= &v_t(t+h)- v_t(t)+\int_0^{t}{\cal R}(\tau)\big[g(u(t+h-\tau))-g(u(t-\tau))\big]d\tau
%\\
%&+\int_t^{t+h}{\cal R}(\tau)g(u(t+h-\tau))d\tau
%\,,
%\end{split}
%\end{equation*}
%and hence \eqref{stima1} and Proposition \ref{eq:lipcon} we obtain
%\begin{equation*}
%\begin{split}
%\|u_t(t+h)- u_t(t)\|_{L^2}
%\le &\|v_t(t+h)- v_t(t)\|_{L^2}+\int_0^{t}\|g(u(s+h))-g(u(s))\|_{L^2}ds
%\\
%&+\int_0^{h}\|g(u(s))\|_{L^2}ds
%\\
%\le &\|v_t(t+h)- v_t(t)\|_{L^2}+C\int_0^{t}\|\nabla u(s+h)-\nabla u(s)\|_{L^2}ds
%\\
%&+C\int_0^{h}\|\nabla u(s)\|_{L^2}ds
%\\
%\le &\|v_t(t+h)- v_t(t)\|_{L^2}+C h
%.
%\end{split}
%\end{equation*}
\begin{equation*}
u_t(t)=v_t(t)+\int_0^t {\cal R}(t-\tau)g(u(\tau))d\tau,
\end{equation*}
and hence, thanks to the regularity of $v(t)$, see \eqref{eq:v}, and $u\in C([0,T];H^1_0(\Omega))$ we get $u\in C^1([0,T];L^2(\Omega))$. Therefore, one can restart by the data $(u(T),u_t(T)\in H^1_0(\Omega)\times L^2(\Omega)$ but this is in contrast with the fact that  $T$ is maximal. The contradiction follows by assuming that $T$ is a positive real number and hence $T=\infty$.

Now  we suppose that the constant $C_0$ in \eqref{eqn:hypnn} satisfies the extra condition \eqref{eq:c0}. By \eqref{eqn:hypnn1} and \eqref{eq:lambda2} we get
\begin{equation}
\int_\Omega G(u)\ dx\le \frac{C_0}{\lambda}\int_\Omega |\nabla u|^2\ dx
\,.
\end{equation}
Therefore, putting the previous estimate into the expression \eqref{eqn:energy} of the energy, we obtain
\begin{equation*}
E(t)\ge\;\frac{1}{2}\int_{\Omega}|u_t|^2\ dx+\Big(\frac{1-a(0)}{2}-\frac{C_0}{\lambda}\Big)\int_{\Omega} |\nabla u|^2\ dx
\,,
\end{equation*}
that is \eqref{eq:coercive}
where $C=\frac{\lambda(1-a(0))-2C_0}{\lambda}>0$ thanks to the assumption $C_0<\lambda(1-a(0))/2$. In particular
$E(0)\ge0$.

Again by \eqref{eq:diss}, we get
\begin{multline}\label{eq:diss1}
E(t)
\le \frac{1}{2}\|u_1\|_{L^2}^2+\frac{1+a(0)}{2}\|\nabla u_0\|_{L^2}^2+\int_\Omega |G(u_0)|\ dx
\\
-a(t)\int_\Omega \nabla u_{0}\cdot \nabla u(t)\ dx
-\int_0^t\dot a(s)\int_\Omega \nabla u_{0}\cdot \nabla u(s)\ dx\ ds
\,.
\end{multline}
If $C>0$ is the constant in \eqref{eq:coercive}, taking into account that
\begin{equation*}
-a(t)\int_\Omega \nabla u_{0}\cdot \nabla u(t)\ dx
\le
\|a\|_\infty\int_\Omega |\nabla u_{0}|\ |\nabla u(t)|\ dx
\le
\frac{C}4\int_{\Omega} |\nabla u|^2\ dx
+\frac{\|a\|_\infty^2}{C}\|\nabla u_0\|_{L^2}^2
\,,
\end{equation*}
\begin{equation*}
-\int_0^t\dot a(s)\int_\Omega \nabla u_{0}\cdot \nabla u(s)\ dx\ ds
\le
\frac{C}4\int_0^t|\dot a(s)|\int_\Omega |\nabla u(s)|^2\ dx\ ds
+\frac{\|\dot a\|_1}{C}\|\nabla u_0\|_{L^2}^2
\,,
\end{equation*}
from  \eqref{eq:diss1}
we get
\begin{multline}\label{eq:diss2}
E(t)
\le \frac{1}{2}\|u_1\|_{L^2}^2+\Big(\frac{1+a(0)}{2}+\frac{\|a\|_\infty^2+\|\dot a\|_1}{C}\Big)\|\nabla u_0\|_{L^2}^2+\int_\Omega |G(u_0)|\ dx
\\
+\frac{C}4\int_{\Omega} |\nabla u|^2\ dx
+\frac{C}4\int_0^t|\dot a(s)|\int_\Omega |\nabla u(s)|^2\ dx\ ds
\,.
\end{multline}
Putting together \eqref{eq:coercive}  and \eqref{eq:diss2}, 
we get
\begin{multline*}
\frac12\| u_t(t)\|_{L^2}^2+\frac{C}4\|\nabla u(t)\|_{L^2}^2
\le \frac{1}{2}\|u_1\|_{L^2}^2+\Big(\frac{1+a(0)}{2}+\frac{\|a\|_\infty^2+\|\dot a\|_1}{C}\Big)\|\nabla u_0\|_{L^2}^2+\int_\Omega |G(u_0)|\ dx
\\
+\frac{C}4\int_0^t|\dot a(s)|\|\nabla u(s)\|_{L^2}^2\ ds
\,.
\end{multline*}
Applying Gronwall lemma, we have for any $t\ge0$
\begin{multline*}
\frac12\| u_t(t)\|_{L^2}^2+\frac{C}4\|\nabla u(t)\|_{L^2}^2
\\
\le 
e^{\|\dot a\|_1}\Big(
\frac{1}{2}\|u_1\|_{L^2}^2+\Big(\frac{1+a(0)}{2}+\frac{\|a\|_\infty^2+\|\dot a\|_1}{C}\Big)\|\nabla u_0\|_{L^2}^2+\int_\Omega |G(u_0)|\ dx\Big)
\,.
\end{multline*}
%This contradicts the maximality of $T$. So,  $T=\infty$. Therefore, $u$ is  global 
Moreover,  putting the above estimate into 
\eqref{eq:diss2} and taking into account \eqref{eq:Gu-0} we obtain that 
\eqref{eq:boundE} holds true.  
Finally, \eqref{eq:boundA}  follows from \eqref{eq:coercive} and \eqref{eq:boundE}, while \eqref{eq:boundG*A} holds for strong solutions in view of
\eqref{eq:disss}. 
\end{Proof}

Under more regular assumptions on the integral kernel, we can establish a different
result concerning the global existence of solutions and the dissipation of energy. 
Indeed, we will assume that the integral kernel satisfies the following conditions 
\begin{equation}\label{eq:k2}
\begin{split}
&a\in C^1 ([0,\infty)),\ \dot{a}(0)<0,\ a(t)\ge0,\ \dot{a}(t)\le0\ \forall t\ge0,
\\
& \ddot{a}(t)\in L^1_{loc}(0,+\infty), \ddot{a}(t)\ge0, \ \text{a.e.}\ t\ge0,
\\
&a\,, \dot{a}\in L^1(0,+\infty),\ a(0)<1\,.
\end{split}
\end{equation} 
It is well known that these conditions imply that $a$ is a strongly positive definite kernel, see \cite[Corollary 2.2]{NS}, and hence \eqref{eq:k1} holds true.
Then we can consider a different expression for the energy of the solutions  with respect to \eqref{eqn:energy}. More precisely, we will define the energy as follows
\begin{equation}\label{eqn:energy1}
\begin{split}
E(t)=&\frac{1}{2}\int_{\Omega} |u_t(t,x)|^2\ dx+\frac{1-a(0)+a(t)}2\int_{\Omega} |\nabla u(t,x)|^2\ dx
\\
&
-\frac{1}{2}\int_{\Omega}\int_0^t \dot a(t-s) |\nabla u(s,x)-\nabla u(t,x)|^2\ ds\ dx
-\int_{\Omega} G(u)\ dx\qquad t\ge0\,.
\end{split}
\end{equation}
Thanks to the assumptions \eqref{eq:k2} $E(t)$ is a decreasing function, see e.g. \cite{MS,ACS}. In particular, we have
\begin{equation}\label{eq:deriv_energy}
E'(t)=
\frac12\dot a(t)\int_{\Omega} |\nabla u(t,x)|^2\ dx
-\frac{1}{2}\int_{\Omega} \int_0^t \ddot a(t-s)|\nabla u(s,x)-\nabla u(t,x)|^2\ ds\ dx
\qquad\ \text{a.e.}\ t\ge0
\,.
\end{equation}

\begin{theorem}\label{th:energy}
Let us assume \eqref{eq:k2}, \eqref{eq:g1}, \eqref{eqn:hypnn} and \eqref{eq:c0}.

For any $u_0\in H^1_0(\Omega)$ and $u_1\in L^2(\Omega)$ there exists a unique mild solution $u$ on $[0,\infty)$ of 
the Cauchy problem
\begin{equation}\label{eq:cauchy}
\begin{cases}
\displaystyle
u_{tt}(t,x) =\triangle u (t,x)+\int_0^t \dot{a} (t-s)\triangle u(s,x)\,ds+g(u(t,x))\,,
\quad t\ge0,\,\, x\in \Omega,
\\
u(t,x)=0\qquad   t\ge0, \,\, x\in\Gamma,
\\
u(0,x)=u_{0}(x),\quad
u_t(0,x)=u_{1}(x),\qquad  x\in \Omega.
\end{cases}
\end{equation}
In addition, if the initial data are more regular, that is $u_0\in H^2(\Omega)\cap H^1_0(\Omega)$ and $u_1\in H^1_0(\Omega)$
the mild solution of \eqref{eq:cauchy} is a strong one.

Moreover, the energy  of the mild solution $u$,
defined by \eqref{eqn:energy1} is positive and  we have for any $t\ge 0$
\begin{eqnarray}
\label{eq:coercive10}
&\displaystyle \frac{1-a(0)}2\|\nabla u(t)\|_{L^2}^2-\int_{\Omega} G(u)\ dx>C\|\nabla u(t)\|_{L^2}^2
\\
\label{eq:coercive1}
& \displaystyle E(t)\ge \;\frac12\| u_t(t)\|_{L^2}^2+C\|\nabla u(t)\|_{L^2}^2\,, &
\\
\label{eq:boundE1}
& 
\displaystyle
E(t)
\le C\Big(\| u_{1}\|_{L^2}^2 +\| \nabla u_{0}\|_{L^2}^2\Big)\,,&
\\
\label{eq:boundA1}
& \| u_t(t)\|_{L^2}^2+\|\nabla u(t)\|_{L^2}^2
\displaystyle
\le C\Big(\| u_{1}\|_{L^2}^2 +\| \nabla u_{0}\|_{L^2}^2\Big)\,,&
\end{eqnarray}
%for any $t\ge 0$, where $ C(R)$
%is a positive, increasing, upper semicontinuous function such that $C(0)=0$.
where the symbol $C$ denotes positive constants, maybe different.

\end{theorem}

\section{Hidden regularity results}
Throughout this section we will assume on the integral kernel and on the nonlinearity the conditions \eqref{eq:k2}, \eqref{eq:g1}, \eqref{eqn:hypnn} and \eqref{eq:c0}.

We will follow the approach pursued in \cite{K, KL2} for linear wave equations without memory and in \cite{LoretiSforza2} for the linear case with memory.
First, we need to introduce a technical lemma, that we will use later. For the sake of completeness we prefer to give all details of the proof, nevertheless  some steps are similar to those of the linear case.
\begin{lemma}\label{le:tech}
Let $u\in H^2_{loc}((0,\infty);H^2(\Omega))$ be a function satisfying the following equation
\begin{equation}\label{eq:stato}
u_{tt}(t,x) =\triangle u (t,x)+\int_0^t \dot{a} (t-s)\triangle u(s,x)\,ds+g(u(t,x))\,,
\quad 
\text{in}
\ \ 
(0,\infty)\times\Omega.
\end{equation}
If $h:\overline{\Omega}\to\R^N$ is a vector field of class $C^1$, then for any fixed $S,T\in\R$, $0\le S<T$, the following identity holds true
\begin{equation}\label{eq:identity}
\begin{split}
&\int_S^T\int_{\Gamma}
\big[2\partial_\nu \big(u+\dot a*u\big)\ h\cdot\nabla \big(u+\dot a*u\big)
-h\cdot\nu |\nabla \big(u+\dot a*u\big)|^2
+h\cdot\nu\ (u_t)^2\big]
\ d\Gamma\ dt
\\
=&2\Big[\int_\Omega u_t\ h\cdot\nabla \big(u+\dot a*u\big)\ dx\Big]_S^T
+\int_S^T\int_\Omega \sum_{j=1}^N \partial_j h_j\ (u_t)^2\ dx\ dt
\\
&-2\int_S^T\int_\Omega u_t\ h\cdot\int_0^t \ddot a(t-s) \big(\nabla u(s)-\nabla u(t)\big) ds\ dx\ dt
-2\int_S^T \dot a(t)\int_\Omega u_t\ h\cdot\nabla u\ dx\ dt
\\
&+2\int_S^T\sum_{i,j=1}^N\int_\Omega
\partial_i  h_j\partial_i \big(u+\dot a*u\big)\partial_j \big(u+\dot a*u\big)\ dx\ dt
-\int_S^T\int_{\Omega}
\sum_{j=1}^N \partial_jh_j\ |\nabla \big(u+\dot a*u\big)|^2\ dx\ dt
\\
&+2\int_S^T\int_{\Omega} g(u(t))h\cdot\nabla\big(u+\dot a*u\big)\ dx\ dt
\,.
\end{split}
\end{equation}
\end{lemma}
\begin{Proof}
To begin with, we multiply the equation \eqref{eq:stato} by 
\begin{equation*}
2h\cdot\nabla\Big( u(t)+\int_0^t\ \dot a(t-s) u(s)\ ds\Big)
%=2\sum_{j=1}^n h_j\partial_j u
\end{equation*}
and integrate over $[S,T]\times\Omega$.
For simplicity, here and in the following we often drop the dependence on the variables.

First, we will handle the term with $u_{tt}$. Indeed,
integrating by parts in the variable $t$ gives
%We set $w=u(t)-\int_0^t\ k(t-s) u(s)\ ds$
\begin{multline}\label{eq:u'}
2\int_S^T\int_\Omega u_{tt}\ h\cdot\nabla \Big( u(t)+\int_0^t\ \dot a(t-s) u(s)\ ds\Big)\ dx\ dt
\\
=2\Big[\int_\Omega u_{t}\ h\cdot\nabla \Big( u(t)+\int_0^t\ \dot a(t-s) u(s)\ ds\Big)\ dx\Big]_S^T
\\
-2\int_S^T\int_\Omega u_{t}\ h\cdot\nabla  u_{t}\ dx\ dt
-2\int_S^T\int_\Omega u_{t}\ h\cdot\nabla\Big(\int_0^t\ \ddot a(t-s) u(s)\ ds+\dot a(0)u(t)\Big)\ dx\ dt
\,.
\end{multline}
Now, we note that, if we integrate by parts in the variable $x$ then we obtain
\begin{equation}\label{eq:u'1}
2\int_\Omega u_{t}\ h\cdot\nabla u_{t}\ dx
=\int_\Omega h\cdot\nabla (u_{t})^2\ dx
=\int_{\Gamma} h\cdot\nu\ (u_{t})^2\ d \Gamma
-\int_\Omega \sum_{j=1}^N \partial_j h_j\ (u_{t})^2\ dx
\,.
\end{equation}
In addition, we can write
\begin{equation}\label{eq:k'star-u}
\begin{split}
\int_0^t\ \ddot a(t-s) u(s)\ ds
=&\int_0^t\ \ddot a(t-s) (u(s)-u(t))\ ds+\int_0^t\ \ddot a(s)u(t)\ ds
\\
=&\int_0^t\ \ddot a(t-s) (u(s)-u(t))\ ds+\dot a(t)u(t)-\dot a(0)u(t)
\,.
\end{split}
\end{equation}
Therefore, plugging \eqref{eq:u'1} and \eqref{eq:k'star-u} into \eqref{eq:u'} yields
\begin{equation}\label{eq:u'2}
\begin{split}
&2\int_S^T\int_\Omega u_{tt}\ h\cdot\nabla \Big( u(t)+\int_0^t\ \dot a(t-s) u(s)\ ds\Big)\ dx\ dt
\\
=&2\Big[\int_\Omega u_{t}\ h\cdot\nabla \Big( u(t)+\int_0^t\ \dot a(t-s) u(s)\ ds\Big)\ dx\Big]_S^T
-\int_S^T\int_{\Gamma} h\cdot\nu\ (u_{t})^2\ d \Gamma\ dt
\\
&+\int_S^T\int_\Omega \sum_{j=1}^N \partial_j h_j\ (u_{t})^2\ dx\ dt
-2\int_S^T\int_\Omega u_{t}\ h\cdot\int_0^t\ \ddot a(t-s) \big(\nabla u(s)-\nabla u(t)\big)\ ds\ dx\ dt
\\
&-2\int_S^T\dot a(t)\int_\Omega u_{t}\ h\cdot\nabla u\ dx\ dt
\,.
\end{split}
\end{equation}
Now, to manage the terms with $\triangle u$, we set
\begin{equation}\label{eq:w}
w(t)=u(t)+\int_0^t\ \dot a(t-s) u(s)\ ds
\,,
\end{equation}
so, we have to evaluate the term
\begin{equation*}
2\int_S^T\int_\Omega
\triangle w\ h\cdot\nabla w\ dx\ dt
\,.
\end{equation*} 
Integrating by parts in the variable $x$ we get
\begin{equation}\label{eq:triangleu}
\begin{split}
&2\int_S^T\int_\Omega
\triangle w\ h\cdot\nabla w\ dx\ dt
\\
=&2\int_S^T\int_{\Gamma}
\partial_\nu w\ h\cdot\nabla w\ d\Gamma\ dt
-2\int_S^T\int_\Omega
\nabla w\cdot\nabla( h\cdot\nabla w)\ dx\ dt
\,.
\end{split}
\end{equation}
We observe that
\begin{equation}\label{eq:triangleu1}
\begin{split}
2\int_\Omega
\nabla w\cdot\nabla( h\cdot\nabla w)\ dx
=&2\sum_{i,j=1}^N\int_\Omega
\partial_i w\ \partial_i ( h_j\partial_j w)\ dx
\\
=&2\sum_{i,j=1}^N\int_\Omega
\partial_i  h_j\partial_i w\partial_j w\ dx
+2\sum_{i,j=1}^N\int_\Omega
h_j\partial_i w\ \partial_j ( \partial_i w)\ dx
\,,
\end{split}
\end{equation}
and
\begin{equation}\label{eq:triangleu2}
\begin{split}
2\sum_{i,j=1}^N\int_\Omega
h_j\partial_i w\ \partial_j ( \partial_i w)\ dx
=&\sum_{j=1}^N\int_\Omega
h_j\ \partial_j ( \sum_{i=1}^N(\partial_i w)^2)\ dx
\\
=&\int_{\Gamma}
h\cdot\nu |\nabla w|^2\ d\Gamma
-\int_{\Omega}
\sum_{j=1}^N \partial_jh_j\ |\nabla w|^2\ dx
\,.
\end{split}
\end{equation}
Therefore, by putting \eqref{eq:triangleu1} and \eqref{eq:triangleu2} into \eqref{eq:triangleu} we obtain
\begin{equation}\label{eq:triangleuF}
\begin{split}
&2\int_S^T\int_\Omega
\triangle w\ h\cdot\nabla w\ dx\ dt
\\
=&2\int_S^T\int_{\Gamma}
\partial_\nu w\ h\cdot\nabla w\ d\Gamma\ dt
-\int_S^T\int_{\Gamma}
h\cdot\nu |\nabla w|^2\ d\Gamma\ dt
\\
&-2\int_S^T\sum_{i,j=1}^N\int_\Omega
\partial_i  h_j\partial_i w\partial_j w\ dx\ dt
+\int_S^T\int_{\Omega}
\sum_{j=1}^N \partial_jh_j\ |\nabla w|^2\ dx\ dt
\,.
\end{split}
\end{equation}
Finally, by  \eqref{eq:u'2} and \eqref{eq:triangleuF}, taking into account \eqref{eq:w} we have the identity \eqref{eq:identity}.
\end{Proof}

\begin{theorem}
Let  $u_0\in H^2(\Omega)\cap H^1_0(\Omega)$, $u_1\in H^1_0(\Omega)$ and $u$ the strong solution of \begin{equation}\label{eq:cauchy1}
\begin{cases}
\displaystyle
u_{tt}(t,x) =\triangle u (t,x)+\int_0^t \dot{a} (t-s)\triangle u(s,x)\,ds+g(u(t,x))\,,
\quad t\ge0,\,\, x\in \Omega,
\\
u(t,x)=0\qquad   t\ge0, \,\, x\in\Gamma,
\\
u(0,x)=u_{0}(x),\quad
u_t(0,x)=u_{1}(x),\qquad  x\in \Omega.
\end{cases}
\end{equation}
If $T>0$, there is a constant $c_0>0$ independent of $T$ such that $u$ satisfies the inequality
\begin{equation}\label{eq:hidden-k0}
\int_0^T\int_{\Gamma} \Big|\partial_\nu u+\dot a*\partial_\nu u\Big|^2d\Gamma dt
\le c_0\int_0^TE(t)\ dt +c_0E(0)
\,,
\end{equation}
where $E(t)$ is the energy of the solution given by \eqref{eqn:energy1}.

Moreover, for a positive constant $c_0=c_0(T)$ we have
\begin{equation}\label{eq:hidden-k}
\int_0^T\int_{\Gamma} \Big|\partial_\nu u+\dot a*\partial_\nu u\Big|^2d\Gamma dt
\le c_0(\|\nabla u_0\|^2_{L^2}+\|u_1\|^2_{L^2})
\,.
\end{equation}

\end{theorem}

\begin{Proof}
To begin with we consider a vector field $h\in C^1(\overline{\Omega};\R^N)$ such that
\begin{equation}\label{eq:h}
h=\nu
\qquad
\text{on}\quad\Gamma\,,
\end{equation}
see e.g. \cite {K} for the construction of such vector field. From now on, we will denote with $c$ positive constants, maybe different. In particular, we have
\begin{equation}\label{eq:h1}
|h(x)|\le c
\quad\text{and}\quad
\sum_{i,j=1}^N |\partial_i  h_j(x)|\ dx\le c\,,
\quad\forall x\in\overline{\Omega}.
\end{equation}
We will apply the identity \eqref{eq:identity} with the vector field $h$ satisfying \eqref{eq:h} and with $S=0$.
First, we observe that 
\begin{equation}\label{eq:mmm}
u_t=0\,,
\qquad
\nabla u=(\partial_\nu u)\nu
\quad
\text{on}
\quad (0,T)\times\Gamma\,.
\end{equation}
For a detailed proof of the second identity see e.g. \cite[Lemma 2.1]{MM}.
Therefore, thanks to \eqref{eq:mmm}
the left-hand side of \eqref{eq:identity} with $S=0$ becomes
\begin{equation*}
\int_0^T\int_{\Gamma} \Big|\partial_\nu u+\dot a*\partial_\nu u\Big|^2d\Gamma dt
\,,
\end{equation*}
and hence  \eqref{eq:identity} can be written as
\begin{equation}\label{eq:identity1}
\begin{split}
&\int_0^T\int_{\Gamma} \Big|\partial_\nu u+\dot a*\partial_\nu u\Big|^2d\Gamma dt
\\
=&2\Big[\int_\Omega u_t\ h\cdot\nabla \big(u+\dot a*u\big)\ dx\Big]_0^T
+\int_0^T\int_\Omega \sum_{j=1}^N \partial_j h_j\ (u_t)^2\ dx\ dt
\\
&-2\int_0^T\int_\Omega u_t\ h\cdot\int_0^t \ddot a(t-s) \big(\nabla u(s)-\nabla u(t)\big) ds\ dx\ dt
-2\int_0^T \dot a(t)\int_\Omega u_t\ h\cdot\nabla u\ dx\ dt
\\
&+2\int_0^T\sum_{i,j=1}^N\int_\Omega
\partial_i  h_j\partial_i \big(u+\dot a*u\big)\partial_j \big(u+\dot a*u\big)\ dx\ dt
-\int_0^T\int_{\Omega}
\sum_{j=1}^N \partial_jh_j\ |\nabla \big(u+\dot a*u\big)|^2\ dx\ dt
\\
&+2\int_0^T\int_{\Omega} g(u(t))h\cdot\nabla\big(u+\dot a*u\big)\ dx\ dt
\,.
\end{split}
\end{equation}
To prove \eqref{eq:hidden-k0} we have to estimate every term on the right-hand side of \eqref{eq:identity1}. Indeed,
\begin{equation}\label{eq:0u-k*u1}
\begin{split}
&2\Big[\int_\Omega u_t\ h\cdot\nabla \big(u+\dot a*u\big)\ dx\Big]_0^T
\\
=&
2\int_\Omega u_t(T)\ h\cdot\nabla \big(u+\dot a*u\big)(T)\ dx-2\int_\Omega u_1\ h\cdot\nabla u_0\ dx
\\
\le&
c\int_\Omega |u_t(T)|^2\ dx+c\int_\Omega  |\nabla \big(u+\dot a*u\big)(T)|^2\ dx
+c\int_\Omega |u_1|^2\ dx+c\int_\Omega |\nabla u_0|^2\ dx
\,.
\end{split}
\end{equation}
We proceed to evaluate for all $t\in [0,T]$ the term $\int_\Omega  |\nabla \big(u+\dot a*u\big)(t)|^2\ dx$, because that evaluation will be also useful later.
Since for all $t\in [0,T]$
\begin{equation*}
\nabla u(t)+\dot a*\nabla u(t)
=\Big(1-a(0)+a(t)\Big)\nabla u(t)+\int_0^t\ \dot a(t-s) \big(\nabla u(s)-\nabla u(t)\big)\ ds
\,,
\end{equation*}
we have
\begin{equation*}
|\nabla \big(u+\dot a*u\big)(t)|^2
\le
2\big(1-a(0)+a(t)\big)^2|\nabla u(t)|^2
+2\Big(\int_0^t\ |\dot a(t-s)| \big|\nabla u(s)-\nabla u(t)\big|\ ds\Big)^2
\,.\end{equation*}
In view  of $\dot a(t)\le0$, $ a(t)\ge0$ and $a(0)<1$ we get
\begin{equation*}
\begin{split}
\Big(\int_0^t\ |\dot a(t-s)| \big|\nabla u(s)-\nabla u(t)\big|\ ds\Big)^2
&\le\int_0^t\ |\dot a(s)|\ ds\int_0^t\ |\dot a(t-s)| \big|\nabla u(s)-\nabla u(t)\big|^2\ ds
\\
&\le-\int_0^t\ \dot a(t-s) \big|\nabla u(s)-\nabla u(t)\big|^2\ ds
\,,
\end{split}
\end{equation*}
and hence
\begin{equation*}
|\nabla \big(u+\dot a*u\big)(t)|^2
\le
2\big(1-a(0)+a(t)\big)|\nabla u(t)|^2
-2 \int_0^t\ \dot a(t-s) \big|\nabla u(s)-\nabla u(t)\big|^2\ ds
\,.
\end{equation*}
Therefore, taking into account the formula \eqref{eqn:energy1} for the energy, by \eqref{eq:coercive1} and \eqref{eq:coercive10}, we get
\begin{equation*}
2\big(1-a(0)\big)\int_{\Omega}|\nabla u(t)|^2
\le c E(t)
\,,
\end{equation*}
\begin{equation*}
2\int_{\Omega} \left( a(t)|\nabla u(t)|^2
- \int_0^t\ \dot a(t-s) \big|\nabla u(s)-\nabla u(t)\big|^2\ ds\right)\ dx
\le 4E(t)
\end{equation*}
and hence
\begin{equation}\label{eq:0u-k*u0}
\int_\Omega|\nabla \big(u+\dot a*u\big)(t)|^2\ dx
\le cE(t)\,.
\end{equation}
By putting \eqref{eq:0u-k*u0} with $t=T$ into \eqref{eq:0u-k*u1} and using again \eqref{eqn:energy1}, we obtain
\begin{equation*}
2\Big[\int_\Omega u_t\ h\cdot\nabla \big(u+\dot a*u\big)\ dx\Big]_0^T
\le
c E(T)+c E(0)
\,,
\end{equation*}
and hence, since the energy $E(t)$ is decreasing, see  \eqref{eq:deriv_energy}, we have
\begin{equation*}%\label{eq:0u-k*u2}
2\Big[\int_\Omega u_t\ h\cdot\nabla \big(u+\dot a*u\big)\ dx\Big]_0^T
\le
c E(0)
\,.
\end{equation*}
Now, we estimate the second term on the right-hand side of \eqref{eq:identity1} by using \eqref{eq:h1}, the expression of energy \eqref{eqn:energy1}  and \eqref{eq:coercive1}, that is
\begin{equation*}%\label{eq:0u-k*u3}
\int_0^T\int_\Omega \sum_{j=1}^N |\partial_j h_j|\ |u_t|^2\ dx\ dt
\le
c\int_0^TE(t)\ dt
\,.
\end{equation*}
In order to bound the term
\begin{equation*}
2\int_0^T\int_\Omega\big| u_t\ h\cdot\int_0^t\ \ddot a(t-s) \big(\nabla u(s)-\nabla u(t)\big)\ ds\big| dx\  dt
\end{equation*}
we note that, thanks also to \eqref{eq:h1}, we have
\begin{multline}\label{eq:firstk'0}
2c\int_0^T\int_\Omega \big|u_t\big|\ \Big|\int_0^t\ \ddot a(t-s) \big(\nabla u(s)-\nabla u(t)\big)ds\Big|  \ dx\ dt
\\
\le
c\int_0^T\int_\Omega |u_t|^2\ dx\ dt
+c\int_0^T\int_\Omega \Big|\int_0^t\ \ddot a(t-s) \big(\nabla u(s)-\nabla u(t)\big)ds\Big|^2  \ dx\ dt
\,.
\end{multline}
To evaluate the second term on the right-hand side of the previous formula, we observe
\begin{equation*}
\begin{split}
\Big|\int_0^t\ \ddot a(t-s) \big(\nabla u(s)-\nabla u(t)\big)ds\Big|^2  
&\le
\Big(\int_0^t\ |\ddot a(t-s)|^{1/2} |\ddot a(t-s)|^{1/2}\big|\nabla u(s)-\nabla u(t)\big|ds\Big)^2  
\\
&\le
\int_0^t\ \ddot a(s)\ ds\int_0^t\  \ddot a(t-s)\big|\nabla u(s)-\nabla u(t)\big|^2 ds
\\
&=
(\dot a(t)-\dot a(0))\int_0^t\  \ddot a(t-s)\big|\nabla u(s)-\nabla u(t)\big|^2 ds
\,.
\end{split}
\end{equation*}
Therefore, in view of $\dot a\le0$ and formula \eqref{eq:deriv_energy}, giving the derivative of the energy, from the above inequality we obtain
\begin{equation}\label{eq:firstk'1}
\begin{split}
&\int_0^T\int_\Omega\Big|\int_0^t\ \ddot a(t-s) \big(\nabla u(s)-\nabla u(t)\big)ds\Big|^2  dx\ dt
\\
\le &-\dot a(0)\int_0^T\int_\Omega\int_0^t\  \ddot a(t-s)\big|\nabla u(s)-\nabla u(t)\big|^2\ ds\ dx\ dt
\le 2\dot a(0)\int_0^T E'(t) dt
\le -2\dot a(0) E(0)
\,.
\end{split}
\end{equation}
Plugging \eqref{eq:firstk'1} into \eqref{eq:firstk'0} and using  \eqref{eq:coercive1} yield
\begin{equation*}%\label{eq:firstk'}
\begin{split}
&2\int_0^T\int_\Omega\big| u_t\ h\cdot\int_0^t\ \ddot a(t-s) \big(\nabla u(s)-\nabla u(t)\big)\ ds\big| dx\ dt
\\
\le
&c\int_0^T\int_\Omega |u_t|^2\ dx\ dt +cE(0)
\le
c\int_0^TE(t)\ dt +cE(0)
\,.
\end{split}
\end{equation*}
Keeping in mind that $\dot a(t)\ge \dot a(0)$ and by using again \eqref{eq:h1} and \eqref{eq:coercive1}, we get
\begin{equation*}
\begin{split}
&-2\int_0^T\dot a(t)\int_\Omega |u_t\ h\cdot\nabla u|\ dx\ dt
\\
\le&
-2\dot a(0)\ c\int_0^T\int_\Omega |u_t| |\nabla u|\ dx\ dt
\le -\dot a(0)\ c\int_0^T\int_\Omega |u_t|^2+ |\nabla u|^2\ dx\ dt
\le c\int_0^T E(t)\ dt
\,.
\end{split}
\end{equation*}
To evaluate the next two terms on the right-hand side of \eqref{eq:identity1} we will use the estimate \eqref{eq:0u-k*u0}. Indeed, as regards the first one, by means of \eqref{eq:h1} we have that
\begin{equation*}
\begin{split}
&\int_0^T\sum_{i,j=1}^N\int_\Omega
|\partial_i  h_j\partial_i \big(u+\dot a*u\big)\partial_j \big(u+\dot a*u\big)|\ dx\ dt
\\
\le &c
\int_0^T\int_\Omega
\Big(\sum_{i=1}^N|\partial_i \big(u+\dot a*u\big)|\Big)^2\ dx\ dt
\le 2^{N-1}c\ \int_0^T\int_{\Omega}\ |\nabla \big(u+\dot a*u\big)|^2\ dx\ dt
\,.
\end{split}
\end{equation*}
Since, from \eqref{eq:0u-k*u0}, we obtain
\begin{equation}\label{eq:1u-k*u0}
\int_0^T\int_\Omega|\nabla \big(u+\dot a*u\big)|^2\ dx\ dt
\le c\int_0^TE(t)\ dt
,
\end{equation}
thus it follows
\begin{equation*}
\int_0^T\sum_{i,j=1}^N\int_\Omega
|\partial_i  h_j\partial_i \big(u+\dot a*u\big)\partial_j \big(u+\dot a*u\big)|\ dx\ dt
\le c\int_0^TE(t)\ dt
\,.
\end{equation*}
In a similar way, thanks again to \eqref{eq:h1} and \eqref{eq:1u-k*u0} we have
\begin{equation*}%\label{eq:u-k*u}
\int_0^T\int_{\Omega}
\sum_{j=1}^N |\partial_jh_j|\ |\nabla \big(u+\dot a*u\big)|^2\ dx\ dt
\le c\int_0^T\int_{\Omega}\ |\nabla \big(u+\dot a*u\big)|^2\ dx\ dt
\le c\int_0^TE(t)\ dt
\,.
\end{equation*}
Finally, to estimate the last term on the right-hand side of \eqref{eq:identity1} first we  use \eqref{eq:h1}
\begin{multline*}
2\int_0^T\int_{\Omega} g(u(t))h\cdot\nabla\big(u+\dot a*u\big)\ dx\ dt
\le c\int_0^T\int_{\Omega} |g(u(t))|^2\ dx\ dt
+c\int_0^T\int_{\Omega} |\nabla\big(u+\dot a*u\big)|^2\ dx\ dt
\,.
\end{multline*}
Since by \eqref{eq:lipcon1} and \eqref{eq:coercive1} we have
\begin{equation*}
\int_0^T\int_{\Omega} |g(u(t))|^2\ dx\ dt
\le c\int_0^T\int_{\Omega} |\nabla u(t)|^2\ dx\ dt
\le c\int_0^TE(t)\ dt
\,,
\end{equation*}
thanks also to \eqref{eq:1u-k*u0}, we obtain
\begin{equation*}
2\int_0^T\int_{\Omega} g(u(t))h\cdot\nabla\big(u+\dot a*u\big)\ dx\ dt
\le
c\int_0^TE(t)\ dt
\,.
\end{equation*}
In conclusion, the previous argumentations show that the sum of all terms on the right-hand side of \eqref{eq:identity1} can be majorized by $c_0\int_0^TE(t)\ dt +c_0E(0)$, with $c_0>0$ independent of $T$, and hence \eqref{eq:hidden-k0} holds true.
In addition, since $E(t)$ is a decreasing function and
$$E(0)=\frac12\|\nabla u_0\|^2_{L^2}+\frac12\|u_1\|^2_{L^2}-\int_\Omega G(u_0)\ dx\,,$$ 
thanks also to \eqref{eq:Gu-0},
\eqref{eq:hidden-k} follows from \eqref{eq:hidden-k0}. 
\end{Proof}

\begin{corollary}\label{co:density}
For any $T>0$ there exists a unique continuous linear map
\begin{equation*}
{\cal L}:H^1_0(\Omega)\times L^2(\Omega)\to L^2((0,T);L^2(\Gamma))
\end{equation*}
such that for any $u_0\in H^2(\Omega)\cap H^1_0(\Omega)$ and $u_1\in H^1_0(\Omega)$, called $u$ the strong solution of \eqref{eq:cauchy1}, we have
\begin{equation*}
{\cal L}(u_0,u_1)=\partial_\nu u
\,.
\end{equation*}

\end{corollary}

\begin{Proof} 
For $u_0\in H^2(\Omega)\cap H^1_0(\Omega)$ and $u_1\in H^1_0(\Omega)$, if we denote by $u$ the strong solution of problem 
\eqref{eq:cauchy1} and apply Lemma \ref{le:unicita} with $X=L^2(\Gamma)$, then for any $T>0$, thanks to 
\eqref{eq:hidden-k} and \eqref{eq:unicita} there exists a constant $c_0=c_0(T,\|\dot a\|_{L^1})>0$ such that 
\begin{equation*}%\label{eq:hidden-k1}
\int_0^T\int_{\Gamma} |\partial_\nu u|^2d\Gamma dt
\le c_0(\|\nabla u_0\|^2_{L^2}+\|u_1\|^2_{L^2})
\,.
\end{equation*}
By density our claim follows.
\end{Proof}

\begin{remark}
For  the mild solution $u$ of \eqref{eq:cauchy1}
%$u_0\in H^1_0(\Omega)$, 
%$u_1\in L^2(\Omega)$ and  
we can introduce the notation $\partial_\nu u$ instead of ${\cal L}(u_0,u_1)$, thanks to
Corollary {\rm \ref{co:density}}. So, for any $T>0$ we have the following trace theorem:
\begin{equation*}
(u_0,u_1)\in H^1_0(\Omega)\times L^2(\Omega)
\Rightarrow
\partial_\nu u\in L^2((0,T);L^2(\Gamma))\,,
\end{equation*}
and there is a positive constant $c_0$ depending on $T$ and $\|\dot a\|_{L^1}$ such that
\begin{equation}\label{eq:hidden-k1}
\int_0^T\int_{\Gamma} |\partial_\nu u|^2d\Gamma dt
\le c_0(\|\nabla u_0\|^2_{L^2}+\|u_1\|^2_{L^2})
\qquad\forall(u_0,u_1)\in H^1_0(\Omega)\times L^2(\Omega)
\,.
\end{equation}
This result does not follow from the usual trace theorems of the Sobolev spaces. For this reason it is called a hidden regularity result. The corresponding inequality \eqref{eq:hidden-k1}
is often called a direct inequality.
\end{remark}

\begin{theorem}
Assume there exists $c_0>0$ independent of  $t$ such that
\begin{equation}\label{eq:pol} 
\int_0^t E(s)\ ds\leq c_0 E(0)\qquad\forall t \geq 0.
\end{equation}
Then, a constant $C>0$ exists  such that for any $u_0\in H^1_0(\Omega)$ and 
$u_1\in L^2(\Omega)$  the mild solution $u$ of \eqref{eq:cauchy1} satisfies
\begin{equation}\label{eq:hidden-k2}
\int_0^\infty\int_{\Gamma} |\partial_\nu u|^2d\Gamma dt
\le C(\|\nabla u_0\|^2_{L^2}+\|u_1\|^2_{L^2})
\,,
\end{equation}
that is
\begin{equation}
\partial_\nu u\in L^2(0,\infty;L^2(\Gamma))
\,.
\end{equation}
\end{theorem}
\begin{Proof} 
In view of  \eqref{eq:hidden-k0} and \eqref{eq:pol}  we have
\begin{equation*}
\int_0^T\int_{\Gamma} \Big|\partial_\nu u+\dot a*\partial_\nu u\Big|^2d\Gamma dt
\le C E(0)
\qquad
\forall T>0
\,,
\end{equation*}
where the constant $C$ is independent of $T$, and hence
\begin{equation*}
\int_0^\infty\int_{\Gamma} \Big|\partial_\nu u+\dot a*\partial_\nu u\Big|^2d\Gamma dt
\le C E(0)
\,.
\end{equation*}
Finally, thanks to \eqref{eq:unicita2} and \eqref{eq:Gu-0} the estimate \eqref{eq:hidden-k2} follows. 
\end{Proof} 
\begin{remark} 
For example, the assumption \eqref{eq:pol} holds if the energy decays exponentially.
Indeed, if 
%The same conclusion (see \eqref{eq:hidden-k2}) holds if instead of  \eqref{eq:Eexp} we assume the following condition
%\begin{equation}\label{eq:pol} 
%\int_0^t E(s)\ ds\leq c_0 E(0)\qquad\forall t \geq 0,
%\end{equation}
%with $c_0$ independent of  $t$.
there exists $m>0$ such that
\begin{equation*}
-\ddot a(t)\le m\ \dot a(t)
\qquad\mbox{for any}\,\, t\ge 0,
\end{equation*}
that is the kernel $-\dot a$ decays exponentially, we can apply
\cite[Theorem 3.5]{ACS} to have that the energy of the mild solution also decays exponentially. Therefore, there exist $\alpha>0$ such that \begin{equation*}
E(t)\le e^{1-\alpha t}E(0)
\qquad
\forall t\ge0
\,.
\end{equation*}
Also  in the case the integral kernel decays polinomially then \eqref{eq:pol} holds
(see \cite{ACS}).
\end{remark} 

\begin{remark} 
If one assumes more regularity on the integral kernel $k=-\dot a$, then it is possible to approach the study of the equation
\begin{equation*}
u_{tt} =\triangle u -\int_0^t k (t-s)\triangle u(s,x)\,ds+g(u)\,,
\end{equation*}
by using the so-called MacCamy's trick, see \cite{MacCamy}.
Adapted to our case, the trick consists in setting 
\begin{equation*}
v=u-k*u\,,
\end{equation*}
to obtain
\begin{equation*}
u=v+\rho_k*v,
\end{equation*}
(where $\rho_k$ is the resolvent kernel of $k$ and has the same regularity of $k$), so $v$ is the solution of the equation
\begin{equation}\label{eq:v}
v_{tt}+\rho_k(0) v_t +\ddot\rho_k*v+\dot\rho_k(0) v=\triangle v+g(v+\rho_k*v)\,.
\end{equation}
However, in  \eqref{eq:v} the terms $\ddot\rho_k*v$ and $\dot\rho_k(0) v$ have a meaning only if $k$, and hence $\rho_k$, is more regular than in our case. For example, a class of kernels fitting our assumptions (see \cite[Corollary 2.2]{NS}), but not suitable for applying the MacCamy's trick is given by
\begin{equation*}
k(t)=k_0e^{-\sqrt t}
\end{equation*}
for a suitable $k_0>0$.

\end{remark}

\section*{Acknowledgments} 
The authors would like to thank the anonymous referee for helpful comments improving the presentation of the paper.

\end{document}